\documentclass[12pt]{article}
\usepackage[cp1251]{inputenc}
\usepackage[russian]{babel}
\usepackage{amsmath,amssymb,amsthm,url}
\input epsf

\begin{document}

\centerline{\bf Простое доказательство изопериметрической теоремы для плоскости
Лобачевского
\footnote{Обновляемая версия: arxiv.org/abs/0911.5319.
Эта заметка под ред. А.Б. Скопенкова составлена составлена
из отзывов В.О. Бугаенко и О.В. Шварцмана на работу [A], которые были переданы Е.И. Алексеевой в декабре 2009 с предложением использовать их для улучшения своей работы.
Отзывы были составлены для Московской математической конференции школьников, см. \url{www.mccme.ru/mmks}.
(Хотя рецензирование на ММКШ анонимное, рецензенты любезно согласились на публикацию настоящей заметки с их именами.)}
}




\bigskip
В этой методической заметке приводятся чёткая формулировка и короткое
доказательство основного результата работы [A] Е.И. Алексеевой (см. ниже), а также
проясняется его связь с изопериметрической теоремой для плоскости Лобачевского.
Доказательство по сути не отличается от приведённого в [A].
Однако ввиду красоты и важности результата короткое доказательство,
освобожденное от ненужных деталей, может быть интересно читателю.
 
\smallskip
{\bf Теорема 1. } [A] {\it На плоскости Лобачевского  среди треугольников $ABC$ с заданными
длинами двух сторон $AB$ и $AC$ максимальную площадь имеет тот, у
которого угол $A$ равен сумме углов  $B$ и $C$.}

\smallskip
Изопериметрическая теорема для плоскости Лобачевского утверждает, что
среди фигур данной площади, ограниченных спрямляемыми кривыми, наибольший периметр имеет круг.
Эта теорема и её многомерные аналоги давно и хорошо известны [S].
Классическому рассуждению Я. Штейнера об изопериметрах (см., например, [K, стр. 19---22] или [P,  стр. 30---31]) недостаёт именно теоремы 1, чтобы оно стало доказательством изопериметрической
теоремы для плоскости Лобачевского. Действительно, в рассуждении Штейнера в силу этой теоремы нужно рассматривать  треугольники $ABC$ с равенством $\angle A=\angle B+\angle C$. Это равенство означает, что на отрезке $BC$ имеется точка $D$, для которой $\angle BAD=\angle B$, $\angle CAD=\angle C$, что эквивалентно равенствам $DB=DA=DC$. Таким образом, точка $A$ лежит на окружности с диаметром $BC$, а это и требуется в рассуждении Штейнера.
Поэтому мы не исключаем, что теорема 1 была известна специалистам или любителям
элементарной математики, хотя никаких ссылок найти не удалось.
 
\smallskip
{\it Доказательство теоремы 1.}
Обозначим через $\alpha,\beta,\gamma$ углы треугольника $ABC$.
Воспользуемся моделью Пуанкаре в круге.
Вершину $A$ поместим в центр модели.
Рассмотрим евклидову окружность $\omega$ и евклидову
прямую, содержащие гиперболические прямые $BC$ и $AB$ соответственно.
Они пересекаются в двух точках $B$ и $B'$
(рис.~1).

\smallskip
\centerline{\epsfbox{hyp.2}}
\centerline{Рисунок 1}
\smallskip

Докажем, что {\it  площадь гиперболического треугольника $ABC$ равна
удвоенной величине евклидова угла $AB'C$,} который мы обозначим через $\tau$.
Действительно, угол между хордой $BC$ и окружностью $\omega$ также равен $\tau$, 
поскольку угол между хордой и касательной равен сооветствующему вписанному углу.
Так как сумма углов евклидова треугольника $ABC$ равна $\alpha+\beta+\gamma+2\tau=\pi$, то
$S(ABC)=\pi-(\alpha+\beta+\gamma)=2\tau.$

Таким образом, треугольник $ABC$ имеет максимальную площадь тогда и
только тогда, когда угол $AB'C$ максимален. Поскольку длины сторон
$AB$ и $AC$ фиксированы, а меняется лишь угол между ними, можно
считать фиксированными точки $A$ и $B$; тогда точка $C$ может
перемещаться по окружности $\psi$ с центром $A$. Очевидно, что угол
$AB'C$ максимален, если евклидова прямая $B'C$ касается окружности
$\psi$ (рис.~2).

\smallskip
\centerline{\epsfbox{hyp.3}}
\centerline{Рисунок 2}
\smallskip

Это, в свою очередь означает, что евклидов угол
$ACB'$~--- прямой. Последнее условие равносильно тому, что
$\pi/2=\angle CAB'+\angle CB'A=\alpha+\tau$. Сопоставив это с
выведенной ранее формулой $S(ABC)=2\tau$ и формулой
$S(ABC)=\pi-\alpha-\beta-\gamma$ для площади треугольника, получаем требуемое
$\alpha=\beta+\gamma$.
QED

\bigskip
[A] Е. Алексеева, Гиперболические треугольники максимальной площади с двумя
заданными сторонами, Мат. Просвещение, 14 (2010), 175-183.
См. также
J. I. Alekseeva, Hyperbolic triangles of the maximum area with two fixed
sides,
http://arxiv.org/abs/0911.5319.

[K] Д.А.Крыжановский, Изопериметры. М.: Физматгиз, 1959.

[P] Протасов В.Ю. Максимумы и минимумы в геометрии. М.: МЦНМО, 2012.

[S] Schmidt E., Beweis der isoperimetrischen Eigenschaft der Kugel im hyperbolischen
und sph\"arischen Raum jeder Dimensionenzahl, Math. Zeitschrift, 49 (1943), 1-109.

\end{document}